\documentclass[a4wide,12pt]{article}

\usepackage{amsmath}
\usepackage{amssymb}
\usepackage{amsfonts}
\usepackage{amsthm}
\usepackage{mathtools}
\usepackage{mathrsfs}
\usepackage{bm}
\usepackage{dsfont}
\usepackage{latexsym}

\usepackage{microtype}
\usepackage{setspace}
\usepackage{indentfirst}
\usepackage{textcomp}
\usepackage{geometry}
\usepackage{lineno}

\usepackage{graphicx}
\usepackage{tcolorbox}
\usepackage{epstopdf}
\usepackage{float}
\usepackage{subcaption}
\usepackage{multirow}

\usepackage{xcolor}

\usepackage{enumerate}
\usepackage[shortlabels]{enumitem}

\usepackage{lscape}

\usepackage{calc}
\usepackage{comment}
\usepackage{listings}

\usepackage{tikz}
\usetikzlibrary{plotmarks}
\usepackage{pgfplots}

\usepackage{algorithmic}
\usepackage[algoruled,boxed,lined]{algorithm2e}

\usepackage{cite}

\usepackage[hidelinks]{hyperref}
\hypersetup{
  colorlinks = true,
  urlcolor   = blue,
  linkcolor  = blue,
  citecolor  = blue
}

\usepackage{cleveref}

\usepackage{fancyhdr}
\fancypagestyle{plain}{%
  \fancyhf{}
  \fancyfoot[C]{\thepage}
  
}

\ifpdf
  \DeclareGraphicsExtensions{.eps,.pdf,.png,.jpg}
\else
  \DeclareGraphicsExtensions{.eps}
\fi

\title{Convergence of replicator-mutator particle systems \\ via bounded Wasserstein distance}

\author{Giacomo Borghi\thanks{Department of Mathematics, Technical University of Munich \& Munich Center for Machine Learning, Munich, Germany (\texttt{giacomo.borghi@tum.de}).}
}

\usepackage{amsopn}

\usepackage{mathrsfs}
\usepackage{mathtools}
\usepackage{amssymb}
\usepackage[shortlabels]{enumitem}
\usepackage{bm}
\usepackage{yfonts}
\usepackage{eufrak}

\newcommand{\Rd}{\mathbb{R}^d}
\newcommand{\RR}{\mathbb{R}}
\newcommand{\x}{\mathbf{x}}
\newcommand{\w}{\mathbf{w}}
\newcommand{\unif}{\mathrm{Unif}}
\newcommand{\law}{\mathrm{Law}}
\newcommand{\lip}{{\mathrm{Lip}}}
\newcommand{\ve}{\varepsilon}
\newcommand{\W}{\mathbb{W}}

\newcommand{\R}{\mathbf{R}}
\newcommand{\M}{\mathbf{M}}

\renewcommand{\d}{\mathrm{d}}

\newcommand{\OX}{\overline{X}}

\newcommand{\X}{\mathbf{X}}
\renewcommand{\j}{\mathbf{j}}

\newtheorem{theorem}{Theorem}[section]
\newtheorem{lemma}[theorem]{Lemma}
\newtheorem{proposition}[theorem]{Proposition}

\theoremstyle{definition}
\newtheorem{definition}[theorem]{Definition}

\theoremstyle{remark}

\newtheorem{assumption}[theorem]{Assumption}

\numberwithin{equation}{section}

\ifpdf
\hypersetup{
  pdftitle={Munuscript - Convergence of replicator-mutator particle systems via bounded Wasserstein distance - Borghi 2026},
  pdfauthor={Giacomo Borghi}
}
\fi

\begin{document}

\maketitle

\begin{abstract}

Interacting particle systems with selection and mutation mechanisms arise in biological modeling, evolutionary optimization, and particle filtering methods. In this work we will consider the paradigmatic replicator–mutator PDE, which combines a diffusive component with a source term that depends on an external potential. We derive quantitative error estimates for the distance between the empirical measure of an unweighted particle system and the solution of the PDE. The proof uses a Wasserstein-1 distance with a bounded ground cost, which allows us to treat unbounded potentials on an unbounded state space. The key tools are the optimal coupling of newly resampled particles and, for the time-continuous limit, the use of a discrete-time killing resampling.
\end{abstract}

\bigskip
\textbf{keywords:} Interacting particle systems, birth-death dynamics, replicator-mutator equation, optimal transport, propagation of chaos

\bigskip
\textbf{MSCcodes:} 65C35, 60B10, 49Q22, 35K58, 62M20

\section{Introduction and main result}
\label{sec:intro}

Complex systems where different interacting particles undergo birth and death dynamics appear in several contexts. In biology, they describe the evolution of a population of cells which reproduce and mutate \cite{perthame2007transport}. In Bayesian statistics, they constitute a large class of particle filter algorithms, particularly effective for multi-modal targets, under the name Sequential Monte Carlo \cite{delmoral2006sequential}. In optimization, such dynamics is at the basis of genetic algorithms \cite{holland1992adaptation}, a popular heuristic involving particle exploration to find global minimizers of non-differentiable functions. 

All such systems can be described at two different levels. One is the microscopic level of the finite $N$ particle systems. One is the statistical level of an evolving probability distribution $\rho\in \mathcal{P}(\Rd)$ typically obtained by taking the many particle limit $N\to \infty$. In this paper we are interested in rigorously studying the relation between these two descriptions with quantitative approximation errors.
We focus on maybe the simplest of such models, that is, the Replicator-Mutator (RM) PDE which describes particles undergoing birth or death according to a fitness potential $V:\Rd \to \RR$, and diffusion with strength $\sigma>0$:
\begin{equation}
\label{eq:pde}
\partial_t \rho_t = \frac{\sigma^2}2 \Delta \rho_t - \left( V - \mathbb{E}_{\rho_t}[V] \right)\rho_t\,.
\end{equation}
The constant $\mathbb{E}_{\rho_t}[V] = \int V(y) \rho_t(\d y)$ ensures the total mass remains constant during the evolution. The RM appears in mathematical biology for the description of cell dynamics \cite{perthame2007transport}, in machine learning as the grazing collision limit of collective hyperparameter learning dynamics \cite{borghi2026pbt}. It can also be regarded as the time-continuous limit of a bootstrap particle filter \cite{candy2007bootstrap}. From a variational point of view, while \eqref{eq:pde} is not a gradient flow, it is comprised of two gradients: the Laplacian corresponds to the $L^2$-Wasserstein gradient of the log-entropy, while the source term corresponds to the Fisher--Rao gradient of the external energy $\int V \rho$. Interestingly, the two energies together form the KL divergence to $\pi \propto e^{-V/\sigma}$, up to a constant.

\begin{algorithm}[t]
\caption{Replicator-mutator particle dynamics}
\label{alg}
\begin{algorithmic}
\STATE{\textbf{Input:} Potential function $V\geq0$, population size $N$, initial distribution $\rho_0\in \mathcal{P}(\Rd)$, 
diffusion coefficient $\sigma\geq0$, time discretization $\tau \in (0,1]$, time horizon $T>0$.}
\STATE{Particle initialization $X^i_0\sim \rho_0$ i.i.d., $i = 1, \dots,N$}
\STATE{$n = 0$}
\WHILE{$n\tau \leq T$}
\STATE{Compute weights $w_n^i = e^{-\tau V(X^i_n)}/ \sum_{j=1}^N e^{-\tau V(X^j_n)}$, $i = 1, \dots, N$}
\FORALL{$i = 1, \dots, N$}
\STATE{Sample $\xi_n^i\sim\mathcal{N}(0,I_d)$}
\STATE{with probability $1  - e^{-\tau V(X^i_n)}$:}
\STATE{\hspace{0.8cm}  sample $J^i_n \sim \mathrm{Categorical}(w^1_n,\dots, w^N_n)$}  
\STATE{\hspace{0.8cm}  $X^i_{n+1} = X_n^{J^i_n} + \sigma\sqrt{\tau} \xi_n^i$}
\STATE{else:}
\STATE{\hspace{0.8cm}  $X_{n+1}^i = X_{n}^i + \sigma\sqrt{\tau}\xi_n^i$ }
\ENDFOR
\STATE{$n = n+1$}
\ENDWHILE
\end{algorithmic}
\end{algorithm}

\subsection{Contribution}

We are interested in connecting \eqref{eq:pde} and the simple time-discrete particle discretization described by Algorithm \ref{alg}.  At every step $n$, particles may undergo independent resampling according to their potential value $V$ following a discrete-time killing resampling mechanism \cite{chopin2022resampling}, also known as $\epsilon$-bootstrap interaction \cite{verge2015parallel}. Gaussian increments are added at every step to discretize the Laplacian operator.  Let $\tau \in (0,1]$ be the time-step. We will show that, as the number of particles $N$ increases, the particle distribution converges to the splitting scheme 
made by alternating distribution reweighting and convolution with the heat kernel $K_{s}$
\[
\R_\tau[\rho] :=\frac{e^{-\tau V}\rho}{\rho(e^{-\tau V})}  \qquad \textup{and} \qquad  
\M_\tau[\rho] := K_{\sigma \sqrt{\tau}} \ast \rho
\]
leading to the time discrete update
\begin{equation}
\label{eq:rhon}
\rho_{n+1} = \M_\tau\left[\R_\tau [\rho_n] \right]\,.
\end{equation}

To study the convergence of the particle system we consider a Wasserstein-1 distance with bounded ground cost:
\[
\W^d_1(\mu, \nu) = \inf_{X\sim \mu, Y \sim \nu}\mathbb{E} \left[ d(X,Y)\right]\,, \qquad \textup{with}\quad d(x,y):= |x - y|\wedge 1\,.
\]
Note that, by the Kantorovich--Rubinstein duality \cite{villani2009}, it corresponds to the Bounded Lipschitz norm between measures (see \eqref{eq:dual}).
Since $d(x,y) \leq |x-y|$, one can use the usual Wasserstein concentration inequalities \cite{fournier2015rate} for i.i.d. particles for the $L^1$-Wasserstein distance
\[
\mathbb{E} \left[ \W^d_1 (\rho_0^N, \rho_0)\right] \lesssim C M_q^{1/q}(\rho_0) \, \ve (N,d)
 \]
 with the rate $\ve(N,d)$ given by 
 \[ \ve(N,d) : = 
\begin{cases}
N^{-1/2} & \textup{if} \;\; d = 1 \\
(1+\log N)N^{-1/2} & \textup{if} \;\; d = 2 \\
N^{-1/d}  & \textup{if} \;\; d \geq 3\,,
\end{cases}
\]
when $q>2$.  Our objective is to show that the rate is preserved at later times, and that the time-discretization introduces an error of order $\tau$.

\begin{assumption}
\label{asm:V} Let $V\in C^2(\Rd)$, $V\geq0$, such that there exist $C_V, C'_V>0$
\begin{equation*}
\begin{split}
V(x) &\leq C_V(1+|x|^2)
 \\
|\nabla V(x) & |\leq C_V' + V(x) \\
|\Delta V(x) | & \leq C^{''}_V (1 + |x|^2)
\end{split} \qquad \textup{for all} \quad x\in \Rd\,.
\end{equation*}
\end{assumption}

\begin{theorem}  \label{t:main}
Given a time horizon $T>0$, let $\rho \in \mathcal{C}([0,T], \mathcal{P}(\Rd)$ be a measure solution to \eqref{eq:pde} (see Definition \ref{def:solution}) with finite $q$-th moments, $q\geq4$. Let $\rho_n^N = 1/N\sum_{i}\delta_{X^i_n}$ be the empirical measure of the particle system constructed according to Algorithm \ref{alg} with time step $\tau \in (0,1]$.

If the potential $V$ satisfies \ref{asm:V}, it holds
\[
\sup_{t_n = n\tau \leq T} \mathbb{E}\left[  \W^d_1 (\rho_n^N, \rho(t_n))\right ] \,  \leq C(1 + T^2) e^{C(1 + T^2)} \left( \, \ve(N,d) \, + \, \tau \right)
\]
with the positive constant $C$ depending on $V,q,d$ and on the $q$-th moment of $\rho_0$.
\end{theorem}

The propagation of chaos result is obtained via the coupling method \cite{diez2022review1}, by considering an auxiliary system of $N$ i.i.d. particles $\{\OX_n^i\}_{i=1}^N$, with $\OX^i_n \sim \rho_n$. Interestingly, the dependence on the dimension $d$ already appears in the coupling estimate:
\[
\mathbb{E} \W_1(\rho_n^N, \overline{\rho}_n^N) \lesssim_T \ve(N,d)\,,
\]
(see Proposition \ref{p:particles}) where mean-field type models, instead, typically exhibit an interaction error of order $1/\sqrt{N}$. This is due to the fact that in Algorithm \ref{alg} the interaction requires, at each step, a resampling procedure. This may explain the empirically observed struggle of discretizing particle filters in high dimension due to the weight collapse phenomenon \cite{bickel2008sharp,snyder2008obstacles}. Convergence of the time-discretization is more standard and is obtained via the Feynman--Kac representation of solutions to \eqref{eq:pde}, as in \cite{integration2013delmoral}.

\subsection{Related literature}

Evolution of probabilities involving selection and mutation, either in the time-continuous form \eqref{eq:pde} or the time-discrete one \eqref{eq:rhon}, has been subject to extensive studies  \cite{DelMoral2004}, with the prominent application being Sequential Monte Carlo methods \cite{delmoral2006sequential}. The time-continuous particle system is sometimes called a \textit{Moran particle process} or a \textit{branching-diffusion process}, and the corresponding nonlinear evolution is called a \textit{Fleming Viot process}  \cite{delmoral2000moran}. What we call the external potential $V$ is sometimes referred to as a \textit{fitness function}.
The evolution can be studied via its Feynman--Kac formulation (see \eqref{eq:fk}) and, typically, particle approximation errors of the type
\[
\sup_{n\tau}\mathbb{E}|\langle \phi, \rho^N_n \rangle - \langle \phi, \rho_n \rangle | \lesssim \frac{1}{\sqrt{N}}
\]
are obtained for a class of observables $\phi \in \Phi$; see, for instance, \cite{review2002crisan}. Note that, differently from Theorem \ref{t:main}, the expectation is taken with respect to one single observable. Intuitively, this corresponds to averaging the error among different runs and considering $\mathbb{E}\rho^N_n$ instead of a single instance $\rho^N_n$.

In the more recent work \cite{journel2022convergence}, instead, convergence is shown with respect to the expected Wasserstein distance by considering evolution on the torus $\mathbb{T}^d$ with rate $\ve(N, d)$, as here. The use of Wasserstein distance with bounded ground cost, as $\W_1^d$, has been used in \cite{tough2022fleming} to study particle systems also evolving in a bounded domain, and undergoing resampling when hitting the boundary. In the work \cite{bo2022probabilistic} the authors also consider a Bounded-Lipschitz (BL) distance, with general RM dynamics which also includes \eqref{eq:pde}. The particle system for which propagation of chaos is established, though, is obtained via the Feynman--Kac representation of the evolution, and, practically, corresponds to a weighted particle discretization. In Algorithm \ref{alg}, instead, we consider unweighted particles, which is computationally more desirable to avoid loss of diversity. In \cite{van2009uniform} the author also considers the BL distance to study the bootstrap particle filter, similar to Algorithm \ref{alg}, and obtains non-quantitative convergence of the particle approximation for non-compact domains and possibly unbounded fitness. Recent work \cite{gerber2026safety} also considers particle filters and establishes almost sure, quantitative, convergence for one-dimensional systems with quadratic potential. Long-time high probability bounds for particle filters are also studied in \cite{del2026perturbation} for non-compact state spaces via total variation and  Kolmogorov distance under confinement assumptions on the potential.
The authors of \cite{erny2022well} consider a jump-diffusion process, as in \eqref{eq:pde}, for an unbounded state space, and use the $L^1$-Wasserstein distance $\W_1$. The analysis assumes directly the resampling mechanism to be Lipschitz with respect to $\W_1$, which does not cover, for instance, quadratic $V$.

Our contribution follows a different line of work on propagation of chaos results for collisional models in kinetic theory, \cite{fournier2016nanbu,cortez2016}, where each particle interacts with a randomly sampled partner from the ensemble. While in Algorithm \ref{alg} there is no binary interaction, the technique of optimally coupling the newly sampled particle is central in our proof. More precisely, we follow the strategy developed in \cite{borghi2026nanbu} for studying Nanbu-type Monte Carlo methods and in \cite{borghi2026chaos} for Genetic Algorithms, where first the limit $N \to \infty$ is considered, and only afterward do we go to the time-continuous settings with $\tau\to 0$.

\section{Preliminary lemmas}
\label{sec:lemmas}

Let us set the notation first.  Given a Borel probability measure $\mu \in \mathcal{P}(\Rd)$, we denote by $M_q(\mu) = \int |x|^q \mu(\d x)$ its $q$-th moment. We equip the space with the Wasserstein-1 distance $\W^d_1$ associated with the cost $d(x,y) = |x - y| \wedge 1 = \min(|x - y|,1)$.
By Kantorovich--Rubinstein duality \cite{villani2009}, $\W^d_1$ corresponds to the Bounded Lipschitz norm between measures
\begin{equation}\label{eq:dual}
\W^d_1(\mu, \nu)  = \| \mu - \nu \|_{\mathrm{BL}} = \sup \left \{ \int \phi (\mu - \nu)\, \middle |\, \|\phi\|_{\mathrm{Lip}} \leq 1, \mathrm{osc}(\phi) \leq 1 \right \} \,,
\end{equation}
where $\|\phi\|_\lip = \sup_{x,y} |\phi(x) - \phi(y)|/|x-y||$ and $\mathrm{osc}(\phi) = \sup \phi - \inf \phi$. Also, note that we can restrict the class to the Lipschitz functions with $\|\phi \|_{\infty} \leq 1/2$ since translations $\phi - c$ do not change the value of the integrals appearing in the dual formulation. 

Random variables are assumed to be defined over an underlying probability space $(\mathbb{P}, \mathcal{F}, \Omega)$ and we write $X \sim \mu$ if $\law(X) = \mu$. For a measurable test function $\phi$ and a measure $\mu$ we write the duality pairing as $\mu(\phi) = \int \phi(x) \mu(\d x)$. 

\subsection{Useful bounds}

The nonlinear update $\rho_n$ defined in \eqref{eq:rhon} admits a useful representation via a time discrete Feynman--Kac formula \cite[Chapter~2]{DelMoral2004}. Let $X_0 \sim \rho_0$ and $(B_t)_{t\geq 0}$ be an independent standard Brownian process and define
\[Y_t = 
X_0 + \sigma B_t\,.
\]
The law $\rho_n$ can be equivalently defined as
\begin{equation}
\label{eq:fk}
\rho_n(\phi)
=
\frac{
\mathbb{E} \left[
\phi(Y_{n\tau})
\exp \left(-\tau\sum_{k=0}^{n-1}V(Y_{k\tau})\right)
\right]}
{
\mathbb{E} \left[
\exp \left(-\tau\sum_{k=0}^{n-1}V(Y_{k\tau})\right)
\right]}
\end{equation}
for any bounded test function $\phi$.

\begin{lemma} \label{l:moment}
Let $\rho_0 \in \mathcal{P}_q(\Rd)$ with $q>2$. If $V$ satisfies Assumption \ref{asm:V}, 
then for every $\tau >0$ there exists a constant $C = C(d,q,\sigma,V,\rho_0)>0$ such that
\[
\sup_{n\tau \leq T} M_q(\rho_n) \leq C(1 + \sqrt{T})e^{C(1 + T^2)} = :\overline{M}_T \]
and
\[
\inf_{n\tau \leq T}  \prod_{k=0}^{n-1} Z_\tau(\rho_k) \geq  Ce^{-C(1 + T^2)} = :
e^{-\overline{V}_T}\,.\]
\end{lemma}
\begin{proof}
We use \eqref{eq:fk} with $\phi(x) = |x|^q$ to obtain
\begin{equation}\label{eq:fkMq}
M_q(\rho_n) =
\frac{
\mathbb{E} \left[
|Y_{n\tau}|^q
\exp \left(-\tau\sum_{k=0}^{n-1}V(Y_{k\tau})\right)
\right]}
{
\mathbb{E} \left[
\exp \left(-\tau\sum_{k=0}^{n-1}V(Y_{k\tau})\right)
\right]}
\end{equation}

 Recall that for a Brownian process, it holds that $\mathbb{E}|B_t|^q = t^{q/2} \mathbb{E}|B_1|^q$. Therefore, we have for all $t\in [0,T]$
\[
\mathbb{E}|Y_t|^q \leq C_{q,d} \left( \mathbb{E}|X_0|^q + \sigma^q\mathbb{E}|B_t|^q\right) \leq C_{d,q}\left( M_q(\rho_0) + T^{q/2}\sigma^q \right) \,.
\]
Since the exponential weights $\exp \left(-\tau\sum_{k=0}^{n-1}V(Y_{k\tau})\right) $ take values in $(0,1]$, the numerator in \eqref{eq:fkMq} can be bounded as  
\[
\mathbb{E} \left [
|Y_{n\tau}|^q
\exp \left (-\tau\sum_{k=0}^{n-1}V(Y_{k\tau})\right) \right] \leq C (1 + T^{q/2})\,,
\]
with $C = C(d,q,\sigma, M_q(\rho_0))>0$. By the quadratic upper bound on $V$ in Assumption \ref{asm:V}, therefore we also have
\[
 \mathbb E V(Y_t)\leq C\left(1+ t \right)\, .
\]
For the denominator in \eqref{eq:fkMq}, instead, we apply Jensen's inequality to obtain 
\begin{equation}\label{eq:fk-lower}
 \mathbb E\left[\exp\left(-\tau\sum_{k=0}^{n-1}V(Y_{k\tau})\right)\right]
 \geq \exp\left(-\tau\sum_{k=0}^{n-1}\mathbb E V(Y_{k\tau})\right)
 \geq e^{-CT(1+T)}\,,
\end{equation}
where $C$ is independent of $T$ and $\tau$. Putting the two estimates together in \eqref{eq:fkMq} gives the first desired bound. 

For the second bound we note that, by using the Feynman--Kac formula with $\phi = e^{-\tau V}$,
\[
Z_\tau(\rho_k) = \frac{
\mathbb{E} \left[
e^{-\tau V(Y_{n\tau})}
\exp \left(-\tau\sum_{m=0}^{k-1}V(Y_{m\tau})\right)
\right]}
{
\mathbb{E} \left[
\exp \left(-\tau\sum_{m=0}^{k-1}V(Y_{m\tau})\right)
\right]} = 
\frac{
\mathbb{E} \left[
\exp \left(-\tau\sum_{m=0}^{k}V(Y_{m\tau})\right)
\right]}
{
\mathbb{E} \left[
\exp \left(-\tau\sum_{m=0}^{n-1}V(Y_{m\tau})\right)
\right]} = :\frac{z_k}{z_{k-1}}\,.
\]
Note that $z_0 = 1$ since $e^0 = 1$.
By multiplying this ratio for $k = 1, \dots, n$, the product simplifies to
\[
\prod_{k=0}^{n-1} Z_{\tau}(\rho_{n}) = \prod_{k=0}^{n-1} \frac{z_k}{z_{k-1}} = \frac{z_{n-1}}{z_0} =  
\mathbb E\left[\exp\left(-\tau\sum_{k=0}^{n-1}V(Y_{k\tau})\right)\right]
\]
for which we already derived a lower bound in \eqref{eq:fk-lower}.
\end{proof}

Thanks to the bound on the moments, we can also estimate the error in approximating the partition function from i.i.d. samples. For a bound of order $\tau/\sqrt N$ at a single step, it suffices to assume $q\geq4$, since $V$ grows at most quadratically.
\begin{lemma} 
\label{l:ZN}
Let $\rho\in \mathcal{P}_q(\Rd)$ with $q\geq 4$ and $X^i\sim \rho$, i.i.d. with empirical measure $\overline{\rho}^N = 1/N \sum_{i=1}^N \delta_{X^i}$. If $V$ satisfies Assumption \ref{asm:V}, it holds that
\[
\mathbb{E} \left| Z_\tau(\overline{\rho}^N) - Z_\tau(\rho) \right| \leq C_V \frac{\tau}{\sqrt{N}}\left ( 1 + M_q(\rho)^{2/q}\right)\,.
\]
\end{lemma}
\begin{proof}
Since $\mathbb{E}[e^{-\tau V(X^i)}] = Z_\tau(\rho)$ and the samples are independent, we have
\[
\mathbb{E}|Z_\tau(\overline{\rho}^N) - Z_\tau(\rho)|^2 = \frac1N \mathrm{Var}_\rho(e^{-\tau V})\,.
\]
By Cauchy--Schwarz, then we have
\[
\mathbb{E}|Z_\tau(\overline{\rho}^N) - Z_\tau(\rho)| \leq \frac{ \mathrm{Var}_\rho(e^{-\tau V})^{1/2}}{\sqrt{N}} \leq 
 \frac{
\big(\int (1 - e^{-\tau V(x)})^2 \rho(\d x)\big)^{1/2} }{\sqrt{N}} \leq \frac{\tau}{\sqrt{N}}\left(\int V^2(x) \rho(\d x)\right)^{1/2}
\]
where we also used $1-e^{-u}\leq u $ with $u = \tau V(x) \geq 0$.
Finally, since $V$ satisfies the quadratic upper bound of Assumption \ref{asm:V}, we have
\[
\left(\int V^2(x) \rho(\d x)\right)^{1/2} \leq C_V\left( 1 + M_4(\rho)^{1/2}\right) \leq C_V\left( 1 + M_q(\rho)^{2/q}\right)\,.
\]
\end{proof}

\subsection{Stability of replicator operator}

We collect another important result concerning the Lipschitz continuity of $\R_\tau$ with respect to the bounded Wasserstein distance $\W^d_1$. 
A key aspect of the estimate is its asymmetry with respect to the two input measures, similarly to \cite[Theorem 15]{sprungk2020local}, as only one of them appears in the Lipschitz constant through its partition function $Z_\tau(\nu)$.
\begin{lemma} \label{l:stability}
Let $V$ satisfy Assumption \ref{asm:V} and let $\mu,\nu\in\mathcal P(\Rd)$. It holds that
\[
\W_1^d\bigl(\R_\tau[\mu],\R_\tau[\nu]\bigr)
\leq \frac{ 1+ \tau C'_V}{Z_\tau(\nu)}\, \W_1^d(\mu,\nu)\, .
\]
\end{lemma}
\begin{proof}
We will make use of the dual formulation of $\W^d_1$ as the bounded Lipschitz norm \eqref{eq:dual}. Let $\phi$ be such that $\|\phi\|_\lip \leq 1$ and $\mathrm{osc}(\phi) \leq 1$ and consider 
\[
\hat \phi = e^{-\tau V}\left( \phi - \R_\tau[\mu](\phi) \right )\,.
\]
Recall $\R_\tau[\mu] = e^{-\tau V}\mu / Z_\tau(\mu)$.
By performing some algebraic computations we can show
\begin{align*}
\R_\tau[\mu](\phi)  - \R_\tau[\nu](\phi) &
 = \R_\tau[\mu](\phi) \cdot \frac{Z_\tau(\nu)}{Z_\tau(\nu)}   - \frac{\nu(e^{-\tau V}\phi)}{Z_\tau(\nu)}
 \pm \frac{\mu(e^{-\tau V}\phi)}{Z_\tau(\nu)} \\
 &  = 
 \frac{\nu\left(e^{-\tau V} \R_\tau[\mu](\phi)\right)
 - \nu(e^{-\tau V}\phi)  \pm \mu(e^{-\tau V}\phi)
 }{Z_\tau(\nu)} \\
 &  = 
 \frac{\nu\left(e^{-\tau V} \R_\tau[\mu](\phi)\right)
 - \nu(e^{-\tau V}\phi)  + \mu(e^{-\tau V}\phi) - \mu(e^{-\tau V} \R_\tau[\mu](\phi) )
 }{Z_\tau(\nu)} \\
 & = \frac{ \mu(\hat \phi) - \nu(\hat \phi)}{Z_\tau(\nu)} \,.
\end{align*}

Note that since $V \geq 0$ and translation by $\R_\tau[\mu](\phi)$ does not change the oscillations of the test function, we have $\mathrm{osc}(\hat \phi) \leq 1$. For the Lipschitz constant, we have $|\nabla \phi| \leq 1$ and by Assumption \ref{asm:V} it holds that $|\nabla V|\leq V + C'_V$, which leads to 
\begin{align*}
|\nabla \hat \phi|  &= e^{-\tau V}|\nabla \phi - \tau (\phi - \R_\tau[\mu](\phi))\nabla V| \\
& \leq  e^{-\tau V} (1 + \tau V) + \tau C'_V e^{-\tau V} \leq ( 1+ \tau C'_V)\,,
\end{align*}
where we also used $e^{-u}(1 + u) \leq 1$ with $u = \tau V$.
Therefore, we have that $\hat \phi /(1 + \tau C'_V)$ is a suitable test function and we have $\mu(\hat \phi) - \nu(\hat \phi)\leq (1 + \tau C'_V)\W_1^d(\mu, \nu)$, from which the desired bound follows.
\end{proof}
We note that, if $\nu(V)<\infty$, Jensen's inequality gives $Z_\tau(\nu)\geq e^{-\tau\nu(V)}$ and the constant is of order $1+\mathcal O(\tau)$. Keeping $Z_\tau(\rho_k)$ in the estimate allows us to apply Lemma \ref{l:moment} directly when iterating over $k$.

\subsection{Stability of resampling probabilities}

The last result we collect concerns the probability of a particle undergoing resampling. Specifically, given two particles $x, y \in \Rd$, we ask ourselves what the probability is that they are simultaneously resampled, or not. From Algorithm \ref{alg}, we have that an $x\in \Rd$ particle is resampled with probability 
\[
\beta(z):= 1 - e^{-\tau V(z)}\,.
\]
If we try to maximize the correspondence between particles $x$ and $y$, we have:
\begin{itemize}[itemsep = 0pt, label = $-$]
\item with probability $\min\{\beta(x),\beta(y)\}$: both resampled
\item with probability $|\beta(x)- \beta(y)|$: only one resampled
\item with probability $ 1- \max\{\beta(x),\beta(y)\}$: neither resampled
\end{itemize}

\begin{lemma} 
\label{l:beta}
Let the potential $V$ satisfy Assumption \ref{asm:V}. For any $x,y \in \Rd$, it holds that 
\[
|\beta(x) - \beta(y)| -  \max\{\beta(x), \beta(y)\} d(x ,y) \leq  \tau C'_V
d(x, y)
 \,.
\]
\end{lemma}

\begin{proof} Consider first the case $d(x,y) = 1$. Since $\beta(x), \beta(y) \in [0,1]$, it holds that
\[
|\beta(x) - \beta(y)| -  \max\{\beta(x), \beta(y)\} \leq 0
\]
and the bound is satisfied.
For $d(x,y)=|x-y|<1$, recall from Assumption \ref{asm:V} that $|\nabla V|\leq C'_V+V$ almost everywhere. We have
\[
|\nabla\beta|\leq \tau(C'_V+V)e^{-\tau V}\leq\tau C'_V+\beta\,,
\]
where in the last step we used $u e^{-u}\leq1-e^{-u}$ with $u=\tau V$. Let $m=\max\{\beta(x),\beta(y)\}$. The function $\min\{\beta,m\}$ agrees with $\beta$ at $x$ and $y$ and has gradient bounded by $\tau C'_V+m$ almost everywhere. Its Lipschitz bound therefore gives
\[
|\beta(x) - \beta(y)| \leq |x - y| \tau C'_V + |x - y|  \max\{\beta(x), \beta(y)\}
\]
Given we are considering the case $d(x,y) = |x - y|< 1$, this inequality corresponds to the desired estimate.
\end{proof}

\section{Coupling and limit $N \to \infty$}
\label{sec:N}

\subsection{Two particle systems}

We show propagation of chaos by coupling the particle system $(\X_n)_{n}$ generated via Algorithm \ref{alg} with a nonlinear system $(\overline{\X}_n)_N$ sharing the same random variables. The key coupling is the one of the newly sampled particles, which will depend on i.i.d. continuous random variables $\alpha^i_n \sim \unif[0,N)$.

To this end, 
given a set of positive weights $\w \in \mathrm{int}(\Delta_N)$
we introduce the index map $\j: [0,N) \to \{1, \dots,N\}$
\begin{equation} \label{eq:j}
\j(\w, \alpha):=\min\Big\{\,j\in\{1,\dots,N\}: \alpha<N \sum_{\ell=1}^j w^\ell\,\Big\}\,.
\end{equation}
If $\alpha \sim \unif[0,N)$, note that $\mathbb{P}(\j(\w,\alpha) = j) = w^j$ and so  $\j(\alpha) \sim \mathrm{Categorical}(w^1, \dots, w^N)$ as in Algorithm \ref{alg}. 

Recall a particle $X^i_n$ gets resampled with probability $1 - e^{-\tau V(X^i_n)}$. We make the activation also explicitly dependent on i.i.d. continuous random variables  $\beta^i_n \sim \unif[0,1)$ so that $X^i_n$ is resampled if $\beta^i_n < 1 -  e^{-\tau V(X^i_n)}$.

Starting from $X^i_0 \sim \rho_0$ i.i.d., the update procedure of Algorithm \ref{alg}  can then be rewritten as 
\begin{equation}
X_{n+1}^i = 
\begin{dcases} 
X_n^{\j(\w_n, \alpha^i_n)}  + \sigma \sqrt{\tau} \xi_n^i  & \textup{if}\quad \beta^i_n < 1-  e^{ - \tau V(X^i_n)} \\
X^i_{n} + \sigma\sqrt{\tau}\xi_n^i 
& \textup{else}
 \end{dcases}
\end{equation}
with  the weights given by 
$
\w_n = (w_n^1, \dots, w_n^N)$, $w^i  = e^{-\tau V(X_n^i)} / \sum_{j=1}^N e^{-\tau V(X_n^j)}\,.
$

To construct the nonlinear particle system we will re-use the random variables $X^i_0, \alpha^i_n, \beta^i_n, \xi^i_n$, and, most importantly, optimally couple the newly sampled random variables thanks to the following result, which we recall for completeness.
\begin{lemma}{\cite[Lemma 3.1]{borghi2026chaos}}
\label{l:coupling}
For $p\geq 1$, let $\W^c_p$ be the Wasserstein distance with a lower semi-continuous cost function
$c(x,y)$, and let $f\in\mathcal{P}(\Rd)$ such that $\W_p^c(f, \delta_0)<+\infty$. Consider
$\x=(x^1,\dots,x^N)\in(\Rd)^N$ and weights $\w=(w^1,\dots,w^N)\in\textup{int}(\Delta_N)$
such that $w^i>0$ and $\sum_{i=1}^N w^i=1$. Define the weighted empirical measure as
$
\mu_{\x,\w}:=\sum_{i=1}^N w^i\,\delta_{x^i},
$
and recall from \eqref{eq:j} the index-map $\j(\w,\cdot)$.

There exists a measurable mapping
\begin{align*}
X_{f}^{*}:(\Rd)^N\times \textup{int}(\Delta_N) \times[0,N) &\to \Rd\\
(\x,\w,\alpha) &\mapsto X_f^{*}(\x,\w,\alpha)
\end{align*}
with the following property: if $\alpha\sim\unif[0,N)$, then the pair
$\big(X_f^{*}(\x,\w,\alpha),\,x^{\mathbf{j}(\w, \alpha)}\big)$ is an optimal coupling
between $f$ and $\mu_{\x,\w}$, namely
\[
\mathbb{E}\left[c\Big(X_f^{*}(\x,\w,\alpha), x^{\mathbf{j}(\w, \alpha)}\Big)^p\right]
=  \left(\W_p^c\big(f,\mu_{\x,\w}\big)\right)^p.
\]
\end{lemma}

For any $n \in \mathbb{N}$, let $\X_{\rho_n}^*$ be the map defined above with $f = \R_\tau[\rho_n]$ for the cost $d$ and $p=1$. The nonlinear particle system is iteratively defined as $\overline{\X}_0 = \X_0$ and
\begin{equation}\label{eq:iid}
\OX_{n+1}^i = 
\begin{dcases}
\X_{\rho_n}^* (\X_n, \w_n, \alpha^i_n)   + \sigma \sqrt{\tau} \xi_n^i
  & \textup{if}\quad \beta^i_n < 1-  e^{- \tau V(\OX^i_n)} \\
 \OX^i_{n} + \sigma\sqrt{\tau}\xi_n^i
& \textup{else}\,.
 \end{dcases}
\end{equation}
Note that, by Lemma \ref{l:coupling}, conditionally on the particles at time $n$, we have 
$
\X_{\rho_n}^* (\X_n, \w_n, \alpha^i_n) \sim \R_\tau[\rho_n]
$
and these draws are conditionally independent since the $\alpha^i_n$ are independent. The conditional law of each update then depends only on $\OX^i_n$. If the $\OX^i_n$ are i.i.d.\ with law $\rho_n$, the law after resampling is $e^{-\tau V}\rho_n+(1-Z_\tau(\rho_n))\R_\tau[\rho_n]=\R_\tau[\rho_n]$. Adding the independent Gaussian increments gives $\rho_{n+1}$. Therefore, by induction, we have that 
\[
\OX^i_n \sim \rho_n \qquad \textup{for all}\; i = 1,\dots,N, \quad \textup{i.i.d.} \,,
\]
with $\rho_n$ being the time-discrete update  $\rho_{n+1} = \M_\tau [ \R_\tau [\rho_n]]$, see \eqref{eq:rhon}.

\subsection{Propagation of chaos}

To show convergence as $N \to \infty$, we track the distance $\mathbb{E} d(X^i_n, \OX^i_n)$ for $n\tau\leq T$ and iterate the estimate via a time-discrete Gr\"onwall argument.

\begin{proposition}\label{p:particles} Let $V$ satisfy Assumption \ref{asm:V}, and $\rho_0 \in \mathcal{P}_q(\Rd)$, $q\geq4$. Let $\rho_n^N$ be the empirical measure associated with the algorithm's particle system and 
$\overline{\rho}_n^N$ the one associated with the i.i.d. system \eqref{eq:iid}.
For every fixed $T>0$ and all $n$ such that $n \tau \leq T$ it holds
\[
\mathbb{E} \W_1^d(\rho_n^N, \overline{\rho}_n^N) \leq C(1+T^2)e^{C(1 + T^2)} \ve(N,d)
\]
with $C>0$ depending only on $V,\sigma,d,q$, and the $q$-th moment of $\rho_0$.
\end{proposition}
\begin{proof} Let the particles' resampling probability be given by $\beta(z) = 1 - e^{-\tau V(z)}$. Since the resampling probabilities of the particles are coupled through the same random variables $\beta^i_n\sim\unif[0,1)$, we have that, conditionally on the particles at time $n$, the probabilities of both, only one, and neither being resampled are respectively
\begin{align*}
\mathbb{P}\left( 
 \beta^i_n < \min\{\beta(X^i_n),\beta(\OX_n^i)\} \right)  &=  \min\{\beta(X^i_n),\beta(\OX_n^i)\} \\
\mathbb{P}\left(  \min\{\beta(X^i_n),\beta(\OX_n^i)\} \leq  \beta^i_n \leq  \max\{\beta(X^i_n),\beta(\OX_n^i)\} \right) & =  |\beta(X^i_n) - \beta(\OX_n^i)| \\
\mathbb{P}\left(  \beta^i_n > \max\{\beta(X^i_n),\beta(\OX_n^i)\} \right)  &=  1 -  \max\{\beta(X^i_n),\beta(\OX_n^i)\} 
\end{align*}
The common Gaussian increment does not change $d$. When only one is resampled, we will use that the maximum transport cost to be paid is $\max_{x,\overline{x}} d(x,\overline{x}) = 1$.

By applying Lemma \ref{l:beta} and optimality of the $\alpha$-coupling we have
\begin{align}
\frac 1N\sum_{i=1}^N \mathbb{E} \Big[ &d(X^i_{n+1}, \OX_{n+1}^i)  \Big]
 \leq  \frac 1N\sum_{i=1}^N  \mathbb{E} \left[  \min\{\beta(X^i_n),\beta(\OX_n^i)\}  d(X^{\j(\w_n, \alpha^i_n)}_{n}, \X_{\rho_n}^* (\X_n, \w_n, \alpha^i_n))  \right] \notag 
 \\
 &\qquad 
 +  \frac 1N\sum_{i=1}^N \mathbb{E} \left[ | \beta(X^i_n) - \beta(\OX_n^i)|   
 + \left (1 - \max\{\beta(X^i_n),\beta(\OX_n^i)\} \right) d(X^i_{n}, \OX_{n}^i) \right] 
 \notag
 \\
& \leq \frac 1N\sum_{i=1}^N  \mathbb{E}\left[ \min\{\beta(X^i_n),\beta(\OX_n^i)\} \W_{1}^d \left(\R_\tau[\rho_n^N], \R_\tau[\rho_n] \right) \right]
\notag
 \\
& \qquad   +  \frac 1N\sum_{i=1}^N  \left(\tau C'_V  \mathbb{E} d(X^i_{n}, \OX_{n}^i) + \mathbb{E} d(X^i_{n}, \OX_{n}^i)\right)
\label{eq:gronwall1}
\,.
 \end{align}
Let us focus on the first term. Since $ \min\{\beta(X^i_n),\beta(\OX_n^i)\} \leq \beta(\OX_n^i)$ and $1/N\sum_{i=1}^N \beta(\OX^i_n) = 1 - Z_\tau(\overline{\rho}^N_n)$, we have
\begin{align*}
 \frac 1N\sum_{i=1}^N  \mathbb{E}\Big [ \min\{\beta(X^i_n),\beta(\OX_n^i)\} &\W_{1}^d \left(\R_\tau[\rho_n^N], \R_\tau[\rho_n] \right) \Big]   \leq 
\frac 1N\sum_{i=1}^N  \mathbb{E}\left[ \beta(\OX_n^i) \W_{1}^d \left(\R_\tau[\rho_n^N], \R_\tau[\rho_n] \right) \right] \\
 & =  \mathbb{E}\left[ \left(1 - Z_\tau(\overline{\rho}^N_n)\right) \W_{1}^d \left(\R_\tau[\rho_n^N], \R_\tau[\rho_n] \right) \right] \\
 & \leq (1 - Z_\tau(\rho_n)) \mathbb{E}\W_{1}^d \left(\R_\tau[\rho_n^N], \R_\tau[\rho_n] \right) + \mathbb{E}\left|
 Z_\tau(\rho_n) - Z_\tau(\overline{\rho}_n^N) \right| \\
 &  \leq (1 - Z_\tau(\rho_n)) \mathbb{E}\W_{1}^d \left(\R_\tau[\rho_n^N], \R_\tau[\rho_n] \right) + \frac{C_V\tau}{\sqrt N}\bigl(1+M_q(\rho_n)^{2/q}\bigr)\\
 & \leq (1 - Z_\tau(\rho_n)) \frac{1 + \tau C'_V}{Z_\tau(\rho_n)}\mathbb{E}\W_{1}^d \left(\rho_n^N, \rho_n \right) + \frac{C_V\tau}{\sqrt N}\bigl(1+M_q(\rho_n)^{2/q}\bigr)\,,
\end{align*}
where in the last steps we have used $\W_1^d\leq1$, Lemma \ref{l:ZN} for the difference between the partition functions and Lemma \ref{l:stability} for the stability of the replicator operator $\R_\tau$.

The Wasserstein distance term $\mathbb{E}\W_{1}^d \left(\rho_n^N, \rho_n \right)$ itself, instead, can be bounded by the triangular inequality as
\begin{align*}
\mathbb{E}\W_{1}^d \left(\rho_n^N, \rho_n \right) & 
\leq \mathbb{E}\W_{1}^d \left(\rho_n^N, \overline{\rho}^N_n \right)
+  \mathbb{E}\W_{1}^d \left(\overline{\rho}_n^N, \rho_n \right) \\
& \leq \frac 1N \sum_{i=1}^N \mathbb{E} d(X_n^i, \overline{X}_n^i)  + 
C M_q(\rho_n)^{1/q} \ve(N,d)\,,
\end{align*}
by the concentration inequalities of i.i.d. particles with respect to the Wasserstein-1 distance $\W_1$ \cite{fournier2015rate}. Note that the Monte Carlo rate $1/\sqrt{N}$ is faster than the concentration rate $\ve(N,d)$ (or equal if $d = 1$) and therefore can be absorbed in an $\mathcal{O}(\ve(N,d))$ term.

Continuing the estimate from \eqref{eq:gronwall1} and collecting the two terms in the distance between the coupled particles, we obtain
\begin{multline}  \label{eq:gronwall-step}
 \frac1N\sum_{i=1}^N\mathbb E d(X_{n+1}^i,\OX_{n+1}^i)
 \leq \frac{1+\tau C'_V}{Z_\tau(\rho_n)}
      \frac1N\sum_{i=1}^N\mathbb E d(X_n^i,\OX_n^i) \\
 +C\,\frac{(1+\tau C'_V)(1-Z_\tau(\rho_n))}{Z_\tau(\rho_n)}
      \overline M_T^{1/q}\ve(N,d)
      +\frac{C_V\tau}{\sqrt N}\bigl(1+\overline M_T^{2/q}\bigr)\,,
\end{multline}
where we used the bound $M_q(\rho_n)\leq\overline M_T$ from Lemma \ref{l:moment}.
To iterate the estimate, we add $C\overline M_T^{1/q}\ve(N,d)$ to both sides of \eqref{eq:gronwall-step}. Since
\[
 1+\frac{(1+\tau C'_V)(1-Z_\tau(\rho_n))}{Z_\tau(\rho_n)}
 = \frac{1+\tau C'_V}{Z_\tau(\rho_n)}-\tau C'_V
 \leq \frac{1+\tau C'_V}{Z_\tau(\rho_n)}\,,
\]
we obtain
\begin{multline*}
 \frac1N\sum_{i=1}^N\mathbb E d(X_{n+1}^i,\OX_{n+1}^i)
 +C\overline M_T^{1/q}\ve(N,d) \\
 \leq \frac{1+\tau C'_V}{Z_\tau(\rho_n)}
 \left[\frac1N\sum_{i=1}^N\mathbb E d(X_n^i,\OX_n^i)
       +C\overline M_T^{1/q}\ve(N,d)\right]
 +\frac{C_V\tau}{\sqrt N}\bigl(1+\overline M_T^{2/q}\bigr)\,.
\end{multline*}
Since $Z_\tau(\rho_k)\leq1$, each multiplier is at least one and every partial product appearing in the iteration is bounded by the full product. By Lemma \ref{l:moment}, we have
\[
 \prod_{k=0}^{n-1}\frac{1+\tau C'_V}{Z_\tau(\rho_k)}
 = \frac{(1+\tau C'_V)^n}{\prod_{k=0}^{n-1}Z_\tau(\rho_k)}
 \leq e^{C'_VT+\overline V_T}\,.
\]
Recall that $X_0^i=\OX_0^i$, so the initial particle distance is zero. Iterating the estimate above and using $n\tau\leq T$ therefore gives
\begin{equation*}
 \frac1N\sum_{i=1}^N\mathbb E d(X_n^i,\OX_n^i)
 +C\overline M_T^{1/q}\ve(N,d) 
 \leq e^{C'_VT+\overline V_T}
 \left[C\overline M_T^{1/q}\ve(N,d)
       +\frac{C_V T}{\sqrt N}\bigl(1+\overline M_T^{2/q}\bigr)\right]\,.
\end{equation*}
Subtracting the sampling term from both sides, we obtain
\begin{equation*}
 \frac1N\sum_{i=1}^N\mathbb E d(X_n^i,\OX_n^i)
 \leq C\overline M_T^{1/q}\left(e^{C'_VT+\overline V_T}-1\right)\ve(N,d) \\
 +\frac{C_V T}{\sqrt N}e^{C'_VT+\overline V_T}
 \bigl(1+\overline M_T^{2/q}\bigr)\,.
\end{equation*}

We note that $N^{-1/2}\leq\ve(N,d)$, and so this Monte Carlo rate can be absorbed in the Wasserstein concentration rate $\ve(N,d)$. By recalling the definition of $\overline{M}_T$ and $\overline{V}_T$ from Lemma \ref{l:moment}, we obtain
\begin{align*}
 \frac1N\sum_{i=1}^N\mathbb E d(X_n^i,\OX_n^i) 
 &\leq C\left[(1+\sqrt T)^{1/q}+T+T(1+\sqrt T)^{2/q}\right]
 e^{C(1+T^2)}\ve(N,d) \\
& \leq C(1+T^2)e^{C(1+T^2)}\ve(N,d)\,,
\end{align*}
where the last inequality follows from $q\geq4$. Since $ \mathbb E\W_1^d(\rho_n^N,\overline\rho_n^N)
\leq \frac1N\sum_{i=1}^N\mathbb E d(X_n^i,\OX_n^i)$, we have obtained the desired estimate.
\end{proof}

\section{Limit $\tau \to 0$ and proof of the main Theorem}

We start by stating what we mean by a measure solution to the RM model \eqref{eq:pde}.

\begin{definition}\label{def:solution}
Let $T>0$ and $\rho_0 \in \mathcal{P}(\Rd)$. We say $\rho \in \mathcal{C}([0,T], \mathcal{P}(\Rd))$ is a weak solution to \eqref{eq:pde} with initial datum $\rho_0$ over $[0,T]$ if for all $\phi \in \mathcal{C}^\infty_c(\Rd)$ and $t\in [0,T]$
\[
\rho_t(\phi) = \rho_0(\phi) + \int_0^T\left( \frac{\sigma^2}2\rho_s(\Delta \rho)
 - \rho_s(V\phi) + \rho_s(V) \rho_s(\phi) \right) \d s\,.\] 
\end{definition}

\subsection{Time-discretization error}

To show the RM splitting scheme $\rho_{n+1} = M_\tau [ \R_\tau [ \rho_n]]$ introduces a first order error with respect to the time discretization step $\tau$, we briefly adapt the proof of \cite{integration2013delmoral} to account for unbounded potentials.
\begin{proposition} \label{p:time}
Let $\rho \in \mathcal{C}([0,T], \mathcal{P}(\Rd))$ be a measure solution to \eqref{eq:pde} with $\rho_0\in \mathcal{P}_q(\Rd)$, $q>2$, and $(\rho_n)_{n\in \mathbb{N}}$ be iteratively defined by \eqref{eq:rhon} with $\tau\in (0,1]$. 
If $V$ satisfies Assumption \ref{asm:V}, it holds
\[
\sup_{t_n = n\tau }\W_1^d(\rho_n, \rho(t_n)) \leq C (1 + T^2)e^{C(1 +T^2)} \tau\,,
\]
with the positive constant $C$ depending on $V, q,d$ and the $q$-th moment of $\rho_0$.
\end{proposition}
\begin{proof}
We rely on the Feynman--Kac representation of both the time-continuous solution and its time-discrete approximation; see \eqref{eq:fk} and \cite{DelMoral2004}. Recall $Y_t = X_0 + \sigma B_t$ with $X_0 \sim \rho$ and $B$ an independent Brownian motion. We have for any bounded test function $\phi$
\begin{align*}
\rho_t(\phi)
=
\frac{
\mathbb{E} \left[
\phi(Y_{t})
\exp \left(-\int_0^t V(Y_{s})\d s\right)
\right]}
{
\mathbb{E} \left[
\exp \left(-\int_0^t V(Y_{s})\d s \right)
\right]}\,, \quad 
\rho_n(\phi)
=
\frac{
\mathbb{E} \left[
\phi(Y_{n\tau})
\exp \left(-\tau\sum_{k=0}^{n-1}V(Y_{k\tau})\right)
\right]}
{
\mathbb{E} \left[
\exp \left(-\tau\sum_{k=0}^{n-1}V(Y_{k\tau})\right)
\right]}\,.
\end{align*}

By integrating by parts on each time interval and applying Itô's
formula, we obtain
\begin{align*}
\int_0^{t_n} V(Y_s)\,\d s
 & - \tau \sum_{k=0}^{n-1} V(Y_{t_k})
=
\sum_{k=0}^{n-1}
\int_{t_k}^{t_{k+1}}
\bigl(V(Y_s)-V(Y_{t_k})\bigr)\,\d s
\\
& =
\sum_{k=0}^{n-1}
\int_{t_k}^{t_{k+1}}
(t_{k+1}-s)\,\d[V(Y_s)]
\\
& =
\sigma \sum_{k=0}^{n-1}
\int_{t_k}^{t_{k+1}}
(t_{k+1}-s)\nabla V(Y_s)\cdot\d B_s 
+\frac{\sigma^2}{2}
\sum_{k=0}^{n-1}
\int_{t_k}^{t_{k+1}}
(t_{k+1}-s)\Delta V(Y_s)\,\d s\, .
\end{align*}
To estimate each term, we use the fact that both $|\nabla V|$ and $|\Delta V|$ grow at most quadratically by Assumption \ref{asm:V}. Since $t_{k+1}-s \leq \tau$ for $s\in [t_k, t_{k+1}]$, by using the Burkholder--Davis--Gundy inequality \cite{kallenberg2021}, we have 
\begin{align*}
\mathbb{E} \left| 
\sigma \sum_{k=0}^{n-1}
\int_{t_k}^{t_{k+1}}
(t_{k+1}-s)\nabla V(Y_s)\cdot\d B_s 
\right| & \leq \tau \sigma 
\mathbb{E}\left |
\int_0^{t_n} \nabla V(Y_s)\cdot\d B_s \right| \\
&\leq   C \tau \mathbb{E} \left( \int_0^t |\nabla V(Y_s)|^2 \right)^{1/2} \\
& \leq C \tau \sqrt{t_n} \mathbb{E} \left[1 + \sup_{s\leq t_n} |Y_s|^2\right] \leq  C\tau \sqrt{t_n} \left(1 + M_2(\rho_0) + t_n \right)
\end{align*}
with constants $C$ that depend on $\sigma$ and the growth assumptions on $V$. Note that we had already derived bounds for the moments of $Y_t$ in the proof of Lemma \ref{l:moment}. Similarly, we also have the bound
\begin{align*}
\mathbb{E} \left |
\frac{\sigma^2}{2}
\sum_{k=0}^{n-1}
\int_{t_k}^{t_{k+1}}
(t_{k+1}-s)\Delta V(Y_s)\,\d s \right |& \leq
C \tau \int _0^{t_n} \mathbb{E} \left [1 + |Y_s|^s \right ] \d s 
\\ & \leq C \tau \left( t_n (1 + M_2(\rho_0)) + t_n^2\right)\,.
\end{align*}
Let us take now a test function from the dual formulation \eqref{eq:dual} of $\W_1^d$, that is, $\phi$ with $\| \phi\|_\lip\leq 1$ and $\mathrm{osc}(\phi) \leq 1$. We use a similar trick to that in Lemma \ref{l:stability} and assume, without loss of generality, that $\|\phi - \rho_{t_n}(\phi)\|_\infty\leq1$.
By also applying the inequality $|e^{-a} - e^{-b}|\leq |a - b|$ for positive $a,b$ we have
\begin{align*}
\rho_{t_n}(\phi) - \rho_{n}(\phi) &= \frac{\mathbb{E}\left[ 
(\phi(Y_t) - \rho_{t_n}(\phi))\left(e^{-\tau\sum_{k=0}^{n-1}V(Y_{t_k})} - e^{-\int_0^{t_n} V(Y_s)\d s} \right)
\right]}{\mathbb{E}\left[ e^{-\tau\sum_{k=0}^{n-1}V(Y_{t_k})}\right]} \\
& \leq \frac{\mathbb{E} \left| \int_0^{t_n} V(Y_s)\,\d s
  - \tau \sum_{k=0}^{n-1} V(Y_{t_k}) \right|}{\Pi_{k=0}^N Z_\tau(\rho_n)} \,.
\end{align*}
Now, to bound the numerator, we use the computation we just performed, while to bound the denominator we use Lemma \ref{l:moment}.  Altogether, since $t_n \leq T$, this leads to 
\[
\rho_{t_n}(\phi) - \rho_{n}(\phi) \leq C (1 + T^2)e^{C(1 +T^2)}\tau\,.
\]
The desired estimate is obtained by taking the supremum over $\phi$ and using the dual formulation \eqref{eq:dual}.
\end{proof}

\subsection{Proof of Theorem \ref{t:main}}

\begin{proof}[Proof of Theorem \ref{t:main}] 
Let us recall the different probability measures we have considered: $\rho_n^N$ is the empirical probability measure generated with Algorithm \ref{alg}, $\rho_n$ is the time-discrete update \eqref{eq:rhon}, while $\overline{\rho}_n^N$ is the empirical measure of the particles \eqref{eq:iid} which are i.i.d. with law $\rho_n$. Finally, $\rho(t_n)$ is the solution to \eqref{eq:pde} at time $t_n = n\tau$.

By the triangular inequality, we have
\[
\mathbb{E} \W_1^d\left( \rho^N_n, \rho(t_n) \right) \leq
\mathbb{E} \W_1^d\left( \rho_n^N, \overline{\rho}^N_n \right) +
\W_1^d\left( \overline{\rho}^N_n, \rho_n \right) +
 \W_1^d\left( \rho_n, \rho(t_n) \right)\,.
\]
By Proposition \ref{p:particles}, the first term is of order $(1 + T^2) e^{(1 + T^2)} \ve(N,d)$. Since $\overline{\rho}_n^N$ is made of i.i.d. particles, \cite{fournier2015rate} holds, and so the second error is also of order $(1 + T^2) e^{(1 + T^2)} \ve(N,d)$ given $M_q(\rho_n)^{1/q} \lesssim (1 + T^2) e^{(1 + T^2)}$ by Lemma \ref{l:moment}. Finally, the last term is of order $(1 + T^2) e^{(1 + T^2)} \tau$ by Proposition \ref{p:time}. Altogether, the rate with respect to the number of particles and time discretization is $\mathcal{O}(\ve(N,d) + \tau)$.
\end{proof}

\subsection*{Aknowledgments}

The author was supported by the Advanced Grant NEITALG of the European
Research Council Executive Agency (ERC) (grant agreement No. 101198055), and by the Munich Center for Machine Learning (MCML)

\bibliographystyle{apalike}
\bibliography{bibfile}

\bigskip
\noindent\begin{minipage}{\textwidth}
\centering
\includegraphics[height=1.5cm]{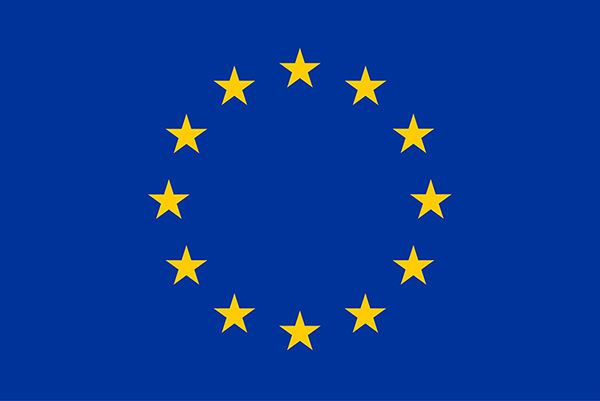}%
\hspace{1em}%
\includegraphics[height=1.5cm]{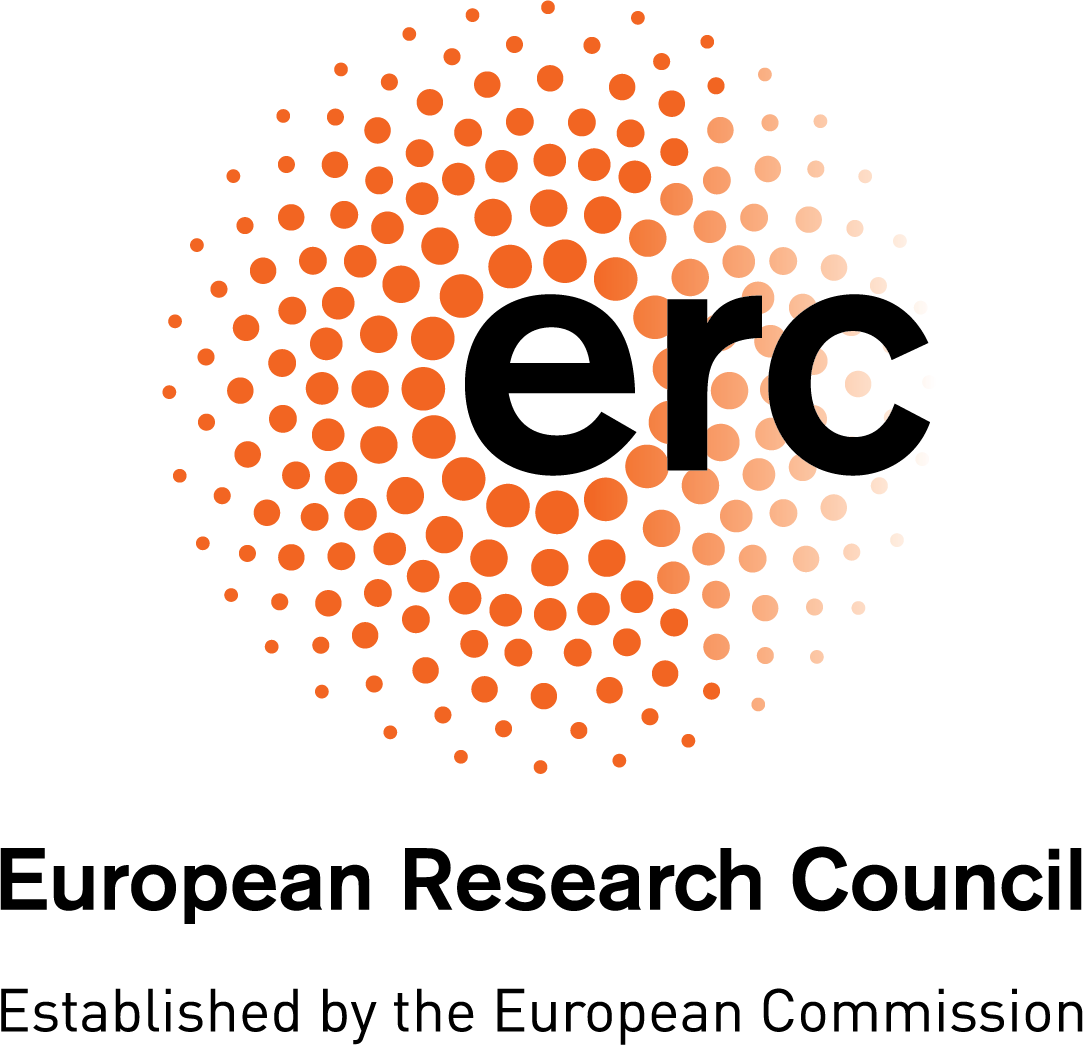}

\vspace{0.5em}
\begin{tcolorbox}[colback=white,boxrule=0.5pt,arc=2pt,left=4pt,right=4pt,top=4pt,bottom=4pt]
\centering\small
Funded by the European Union. Views and opinions expressed are however those of the author(s) only and do not necessarily reflect those of the European Union or the European Research Council Executive Agency. Neither the European Union nor the granting authority can be held responsible for them. This project has received funding from the European Research Council (ERC) under the European Union’s Horizon Europe research and innovation programme (grant agreement No. 101198055, project acronym NEITALG).
\end{tcolorbox}
\end{minipage}

\end{document}